\documentclass[12pt,draft]{article}
\usepackage{amsmath,amscd}
\usepackage{amssymb,latexsym,amsthm}
\usepackage{color}
\usepackage[spanish,english]{babel}
\usepackage{amssymb}
\usepackage{color}
\usepackage{amsmath,amsthm,amscd}
\usepackage[latin1]{inputenc}
\usepackage{lscape}
\usepackage{amsfonts}
\numberwithin{equation}{section}
\newtheorem{theorem}{Theorem}[section]

\newtheorem{proposition}[theorem]{Proposition}
\newtheorem{lemma}[theorem]{Lemma}

\newtheorem{corollary}[theorem]{Corollary}

\theoremstyle{definition}

\newtheorem{definition}[theorem]{Definition}
\newtheorem{examples}[theorem]{Examples}

\newtheorem{remark}[theorem]{Remark}

\newcommand{\cO}{\mbox{${\cal O}$}}

\newcommand{\cU}{\mbox{${\cal U}$}}

\newcommand{\cW}{\mbox{${\cal W}$}}

\title{\textbf{$\sigma$-PBW Extensions of Skew $\Pi$-Armendariz Rings}}
\author{Armando Reyes\footnote{Departamento de  Matem\'aticas. e-mail: mareyesv@unal.edu.co} \\ Universidad Nacional de Colombia, Sede Bogot\'a\\ }
\date{}
\begin{document}
\maketitle
\begin{abstract}
\noindent In this paper we present the notion of skew $\Pi$-Armendariz for the non-commutative rings known as $\sigma$-PBW extensions. This concept  generalizes several definitions of Armendariz rings presented in the literature for these extensions, and in particular, for Ore extensions of injective type. We investigate the relations between skew $\Pi$-Armendariz, $\Sigma$-rigid, $(\Sigma,\Delta)$-skew Armendariz and $\Sigma$-skew Armendariz rings rings.

\bigskip


\bigskip

\bigskip

\end{abstract}
%%%%%%%%%%%%%%%%%%%%%%%%%%%%%%%%%%%%%%%%%%%%%%%%%%%%%%%%%%%%%%%%%%%%%%%%%%%%%
\section{Introduction}\label{section}

In commutative algebra, a ring $B$ is called {\em Armendariz} (the term was introduced by Rege and Chhawchharia in \cite{RegeChhawchharia1997}), if whenever polynomials $f(x)=a_0+a_1x+\dotsb + a_nx^n$, $g(x)=b_0+b_1x+\dotsb + b_mx^m\in B[x]$ such that  $f(x)g(x)=0$, then $a_ib_j=0$, for every $i,j$. The interest of this notion lies in its natural and its useful role in understanding the relation between the annihilators of the ring $B$ and the annihilators of the polynomial ring $B[x]$. As a matter of fact, in  \cite{Armendariz1974}, Lemma 1, Armendariz showed that a  reduced ring (a ring has no nonzero nilpotent elements) always satisfies this condition. In the context of commutative rings and non-commutative rings, more exactly Ore extensions, several treatments have been formulated, see \cite{Armendariz1974}, \cite{RegeChhawchharia1997}, \cite{AndersonCamillo1998},  \cite{KimLee2000}, \cite{Huhetal2002}, \cite{HongKimKwak2003}, \cite{LeeWong2003}, \cite{Matczuk2004} (a detailed list of these works can be found in \cite{NinoReyes2017}, \cite{ReyesSuarez2016b} and  \cite{ReyesSuarezClifford2017}). \\

With the aim of generalizing the results established about Armendariz properties in the mentioned papers above, in this article we are interested in a class of non-commutative rings of polynomial type more general than iterated Ore extensions (of injective type), the $\sigma$-{\em PBW} {\em extensions} (also known as {\em skew Poincar\'e-Birkhoff-Witt extensions}), introduced in \cite{LezamaGallego2011} (see  \cite{Reyes2013PhD} and  \cite{LezamaReyes2014} for a list of non-commutative rings which are $\sigma$-PBW extensions but not iterated Ore extensions). Actually, skew PBW extensions are more general than several fa\-mi\-lies of non-commutative rings such as the following (see \cite{ReyesSuarezClifford2017} for a detailed reference of every one of these families): (i) universal enveloping algebras of finite dimensional Lie algebras; (ii) PBW extensions introduced by Bell and Goodearl; (iii) almost normalizing extensions defined by McConnell and Robson; (iv) sol\-va\-ble polynomial rings introduced by Kandri-Rody and Weispfenning; (v) diffusion algebras studied by Isaev, Pyatov, and Rittenberg; (vi) 3-dimensional skew polynomial algebras introduced by Bell and Smith; (vii) the regular graded algebras studied by Kirkman, Kuzmanovich, and Zhang, and other non-commutative algebras of polynomial type. The importance of skew PBW extensions is that the coefficients do not necessarily commute with the variables, and these coefficients are not necessarily elements of fields (see Definition \ref {gpbwextension} below). In fact, the $\sigma$-PBW extensions contain well-known groups of algebras such as some types of $G$-algebras studied by Levandovskyy and some PBW algebras defined by Bueso et. al., (both $G$-algebras and PBW algebras take coefficients in fields and assume that coefficientes commute with variables),  Auslander-Gorenstein rings, some Calabi-Yau and skew Calabi-Yau algebras, some Artin-Schelter regular algebras, some Koszul algebras, quantum polynomials, some quantum universal enveloping algebras, and others (see \cite{Reyes2013PhD},
\cite{LezamaReyes2014}, \cite{SuarezLezamaReyes2015},  \cite{ReyesSuarezMomento2017},   \cite{SuarezReyesgenerKoszul2017}, or \cite{SuarezReyes2017} for a list of examples). For more details about the relation between $\sigma$-PBW extensions and another algebras with PBW bases, see \cite{Reyes2013PhD}, \cite{LezamaReyes2014}, \cite{ReyesSuarez2017FEJM} and \cite{ReyesSuarezMomento2017}.\\

Since Ore extensions of injective type are particular examples of $\sigma$-PBW extensions (see Example \ref{mentioned}), and ha\-ving in mind that several ring, module and homological properties of these more general rings have been studied by the author and others (see \cite{LezamaGallego2011},  \cite{Reyes2013PhD},   \cite{LezamaAcostaReyes2015},  \cite{Reyes2015},  \cite{ReyesSuarezUMA2017},  \cite{ReyesSuarezClifford2017},  \cite{ReyesSuarezUPTC2016},
\cite{ReyesSuarez2017FEJM}, \cite{ReyesSuarezMomento2017}, \cite{ReyesSuarezClifford2017},  \cite{ReyesSuarezUMA2017}, \cite{SuarezReyesgenerKoszul2017}, \cite{SuarezReyes2017}, \cite{RodriguezReyes2017}, etc), and of course, several notions of Armendariz rings for these extensions have been formulated (see \cite{NinoReyes2017} and \cite{ReyesSuarez2016b}), we consider relevant to formulate a definition of Armendariz ring which extends the proposals developed in these two papers. This is the key objective of our Definition \ref{skewpiArmendariz}, because using this new notion, we show that all before treatments are particular cases of the theory developed here, and hence our results generalize their corresponding assertions in previous papers, including the assertions concerning Ore extensions. It is important to say that our Definition \ref{skewpiArmendariz} is motivated by the theory developed by Luang at. al., in \cite{LunqunJinwangyYueming2013} for the case of Ore extensions. \\

The paper is organized as follows: In Section \ref{definitionexamplesspbw} we establish some useful results about $\sigma$-PBW extensions for the rest of the paper. Section \ref{reunionMelendez} contains the central results of the paper. More exactly, in this section we introduce the notion of {\em skew} $\Pi$-{\em Armendariz ring} (Definition \ref{skewpiArmendariz}) which generalize the notions of $\Sigma$-rigid rings (\cite{Reyes2015}, Definition 3.2), $(\Sigma,\Delta)$-skew Armendariz rings (\cite{NinoReyes2017}, Definition 3.4), $\Sigma$-skew Armendariz rings (\cite{ReyesSuarez2016b}, Definition 3.1), assuming a condition of compatibility introduced in \cite{ReyesSuarezUMA2017}, Definition 3.2. In this section we also study the relation of skew $\Pi$-Armendariz rings and reversible rings (Theorem \ref{hijoamado}). In this way we  generalize some of the results presented by Ouyang in \cite{Ouyang2010} concerning Ore extensions. We conclude with the Remark \ref{masabiert} where we express a possible line of research continuing this work.\\

Throughout the paper, the word ring means a ring (not necessarily commutative) with unity. $\mathbb{C}$ will denote the field of complex numbers, and the letter $\Bbbk$ will denote any field. We recall that a ring $B$ is called {\em reversible}, if $ab=0 \Rightarrow ba=0$, for every $a, b\in B$, and $B$ is called {\em semicommutative}, if $ab=0$ implies $aBb=0$, for every $a, b\in B$ (note that reversible implies semicommutative). For a ring $B$, ${\rm nil}(R)$ represents the set of nilpotent elements of $B$, and $P(B)$ denotes the prime radical of $B$ (the intersection of all prime ideals). $B$ is called a {\em NI} ring, if ${\rm nil}(B)$ forms an ideal of $B$ (from \cite{LiuZhao2006}, Lemma 3.1, we know that if $B$ is a semicommutative ring then the set ${\rm nil}(B)$ is an ideal of $B$).
%Finally, $B$ is called 2-{\em primal}, if $P(B)={\rm nil}(B)$.
%%%%%%%%%%%%%%%%%%%%%%%%%%%%%%%%%%
%%%%%%%%%%%%%%%%%%%%%%%%%%%%%%%%%%
\section{Skew PBW extensions}\label{definitionexamplesspbw}
In this section we establish some useful results about skew PBW extensions for the rest of the paper.
%%%%%%%%%%%%%%%%%
\begin{definition}[\cite{LezamaGallego2011}, Definition 1]\label{gpbwextension}
Let $R$ and $A$ be rings. We say that $A$ is a {\em skew PBW extension} (also known as {\em $\sigma$-PBW extension}) {\em of}  $R$, which is denoted by $A:=\sigma(R)\langle
x_1,\dots,x_n\rangle$, if the following conditions hold:
\begin{enumerate}
\item[\rm (i)]$R\subseteq A$;
\item[\rm (ii)]there exist elements $x_1,\dots ,x_n\in A$ such that $A$ is a left free $R$-module, with basis ${\rm Mon}(A):= \{x^{\alpha}=x_1^{\alpha_1}\cdots
x_n^{\alpha_n}\mid \alpha=(\alpha_1,\dots ,\alpha_n)\in
\mathbb{N}^n\}$,  and $x_1^{0}\dotsb x_n^{0}:=1\in {\rm Mon}(A)$.

\item[\rm (iii)]For each $1\leq i\leq n$ and any $r\in R\ \backslash\ \{0\}$, there exists an element $c_{i,r}\in R\ \backslash\ \{0\}$ such that $x_ir-c_{i,r}x_i\in R$.
\item[\rm (iv)]For any elements $1\leq i,j\leq n$, there exists $c_{i,j}\in R\ \backslash\ \{0\}$ such that $x_jx_i-c_{i,j}x_ix_j\in R+Rx_1+\cdots +Rx_n$.
\end{enumerate}
\end{definition}
%%%%%%%%%%%%%%%%%%%%%%%%%%%%%%%%%%%%%%%%%%%%%%%%
\begin{proposition}[\cite{LezamaGallego2011}, Proposition
3]\label{sigmadefinition}
Let $A$ be a skew PBW  extension of $R$. For each $1\leq i\leq
n$, there exist an injective endomorphism $\sigma_i:R\rightarrow
R$ and an $\sigma_i$-derivation $\delta_i:R\rightarrow R$ such that $x_ir=\sigma_i(r)x_i+\delta_i(r)$, for  each $r\in R$. From now on, we write  $\Sigma:=\{\sigma_1,\dotsc, \sigma_n\}$, and $\Delta:=\{\delta_1,\dotsc, \delta_n\}$.
\end{proposition}
%%%%%%%%%%%%%%%%%%%%%%%%%
%\begin{definition}[\cite{LezamaGallego2011}, Definition 4; \cite{LezamaAcostaReyes2015}, Definition 2.3 (ii)]\label{sigmapbwderivationtype}
%Let $A$ be a skew PBW extension of $R$.
%\begin{enumerate}
%\item[\rm (a)] $A$ is called \textit{quasi-commutative} if the conditions {\rm(}iii{\rm)} and {\rm(}iv{\rm)} in Definition \ref{gpbwextension} are replaced by the following: (iii') for each $1\leq i\leq n$ and all $r\in R\ \backslash\ \{0\}$, there exists $c_{i,r}\in R\ \backslash\ \{0\}$ such that $x_ir=c_{i,r}x_i$; (iv') for any $1\leq i,j\leq n$, there exists $c_{i,j}\in R\ \backslash\ \{0\}$ such that $x_jx_i=c_{i,j}x_ix_j$.
%\item[\rm (b)] $A$ is called \textit{bijective}, if $\sigma_i$ is bijective for each $1\leq i\leq n$, and $c_{i,j}$ is invertible, for any $1\leq i<j\leq n$.
%\item[\rm (c)] $A$ is called of {\em endomorphism type}, if $\delta_i=0$, for every $i$.  In addition, if every $\sigma_i$ is bijective, $A$ is a skew PBW extension of {\em automorphism type}.
%\end{enumerate}
%\end{definition}
%%%%%%%%%%%%%%%%%%%%%%%%%%%%%%%%%%%%%%%%%%%%%%
\begin{examples}\label{mentioned}
If $R[x_1;\sigma_1,\delta_1]\dotsb [x_n;\sigma_n,\delta_n]$ is an iterated Ore extension where
\begin{itemize}
\item $\sigma_i$ is injective, for $1\le i\le n$;
\item $\sigma_i(r)$, $\delta_i(r)\in R$, for every $r\in R$ and $1\le i\le n$;
\item $\sigma_j(x_i)=cx_i+d$, for $i < j$, and $c, d\in R$, where $c$ has a left inverse;
\item $\delta_j(x_i)\in R + Rx_1 + \dotsb + Rx_n$, for $i < j$,
\end{itemize}
then $R[x_1;\sigma_1,\delta_1]\dotsb [x_n;\sigma_n, \delta_n] \cong \sigma(R)\langle x_1,\dotsc, x_n\rangle$ (\cite{LezamaReyes2014}, p. 1212). Note that $\sigma$-PBW extensions of endomorphism type are more general than iterated Ore extensions of the form $R[x_1;\sigma_1]\dotsb [x_n;\sigma_n]$, and in general, $\sigma$-PBW extensions are more general than Ore extensions of injective type (see \cite{LezamaGallego2011},  \cite{Reyes2013PhD} or \cite{LezamaReyes2014} for different examples).
\end{examples}
%%%%%%%%%%%%%%%%%%%%%%%%%%%%%%%%%%%%%%%%%%%%%%
\begin{definition}[\cite{LezamaGallego2011}, Definition 6]\label{definitioncoefficients}
Let $A$ be a skew PBW extension of $R$. Then:
\begin{enumerate}
\item[\rm (i)]for $\alpha=(\alpha_1,\dots,\alpha_n)\in \mathbb{N}^n$,
$\sigma^{\alpha}:=\sigma_1^{\alpha_1}\cdots \sigma_n^{\alpha_n}$,
$|\alpha|:=\alpha_1+\cdots+\alpha_n$. If
$\beta=(\beta_1,\dots,\beta_n)\in \mathbb{N}^n$, then
$\alpha+\beta:=(\alpha_1+\beta_1,\dots,\alpha_n+\beta_n)$.
\item[\rm (ii)]For $X=x^{\alpha}\in {\rm Mon}(A)$,
$\exp(X):=\alpha$, $\deg(X):=|\alpha|$, and $X_0:=1$. The symbol $\succeq$ will denote a total order defined on ${\rm Mon}(A)$ (a total order on $\mathbb{N}^n$). For an
 element $x^{\alpha}\in {\rm Mon}(A)$, ${\rm exp}(x^{\alpha}):=\alpha\in \mathbb{N}^n$.  If
$x^{\alpha}\succeq x^{\beta}$ but $x^{\alpha}\neq x^{\beta}$, we
write $x^{\alpha}\succ x^{\beta}$. Every element $f\in A$ can be expressed uniquely as $f=a_0 + a_1X_1+\dotsb +a_mX_m$, with $a_i\in R$, and $X_m\succ \dotsb \succ X_1$ (eventually, we will use expressions as $f=a_0 + a_1Y_1+\dotsb +a_mY_m$, with $a_i\in R$, and $Y_m\succ \dotsb \succ Y_1$). With this notation, we define ${\rm
lm}(f):=X_m$, the \textit{leading monomial} of $f$; ${\rm
lc}(f):=a_m$, the \textit{leading coefficient} of $f$; ${\rm
lt}(f):=a_mX_m$, the \textit{leading term} of $f$; ${\rm exp}(f):={\rm exp}(X_m)$, the \textit{order} of $f$; and
 $E(f):=\{{\rm exp}(X_i)\mid 1\le i\le t\}$. Note that $\deg(f):={\rm max}\{\deg(X_i)\}_{i=1}^t$. Finally, if $f=0$, then
${\rm lm}(0):=0$, ${\rm lc}(0):=0$, ${\rm lt}(0):=0$. We also
consider $X\succ 0$ for any $X\in {\rm Mon}(A)$. For a detailed description of monomial orders in skew PBW  extensions, see \cite{LezamaGallego2011}, Section 3.
\end{enumerate}
\end{definition}
%%%%%%%%%%%%%%%%%%%%%%%%%%%%%%%%%%%%%%%%%
\begin{proposition}[\cite{LezamaGallego2011}, Theorem 7]\label{coefficientes}
If $A$ is a polynomial ring with coefficients in $R$ with respect to the set of indeterminates $\{x_1,\dots,x_n\}$, then $A$ is a skew PBW  extension of $R$ if and only if the following conditions hold:
\begin{enumerate}
\item[\rm (i)]for each $x^{\alpha}\in {\rm Mon}(A)$ and every $0\neq r\in R$, there exist unique elements $r_{\alpha}:=\sigma^{\alpha}(r)\in R\ \backslash\ \{0\}$, $p_{\alpha ,r}\in A$, such that $x^{\alpha}r=r_{\alpha}x^{\alpha}+p_{\alpha, r}$,  where $p_{\alpha ,r}=0$, or $\deg(p_{\alpha ,r})<|\alpha|$ if
$p_{\alpha , r}\neq 0$. If $r$ is left invertible,  so is $r_\alpha$.
\item[\rm (ii)]For each $x^{\alpha},x^{\beta}\in {\rm Mon}(A)$,  there exist unique elements $c_{\alpha,\beta}\in R$ and $p_{\alpha,\beta}\in A$ such that $x^{\alpha}x^{\beta}=c_{\alpha,\beta}x^{\alpha+\beta}+p_{\alpha,\beta}$, where $c_{\alpha,\beta}$ is left invertible, $p_{\alpha,\beta}=0$, or $\deg(p_{\alpha,\beta})<|\alpha+\beta|$ if
$p_{\alpha,\beta}\neq 0$.
\end{enumerate}
\end{proposition}
%%%%%%%%%%%%%%%%%%%%%%%%
\begin{proposition}\label{lindass}
If $\alpha=(\alpha_1,\dotsc, \alpha_n)\in \mathbb{N}^{n}$ and $r$ is an element of $R$, then
\begin{align*}
x^{\alpha}r = &\ x_1^{\alpha_1}x_2^{\alpha_2}\dotsb x_{n-1}^{\alpha_{n-1}}x_n^{\alpha_n}r = x_1^{\alpha_1}\dotsb x_{n-1}^{\alpha_{n-1}}\biggl(\sum_{j=1}^{\alpha_n}x_n^{\alpha_{n}-j}\delta_n(\sigma_n^{j-1}(r))x_n^{j-1}\biggr)\\
+ &\ x_1^{\alpha_1}\dotsb x_{n-2}^{\alpha_{n-2}}\biggl(\sum_{j=1}^{\alpha_{n-1}}x_{n-1}^{\alpha_{n-1}-j}\delta_{n-1}(\sigma_{n-1}^{j-1}(\sigma_n^{\alpha_n}(r)))x_{n-1}^{j-1}\biggr)x_n^{\alpha_n}\\
+ &\ x_1^{\alpha_1}\dotsb x_{n-3}^{\alpha_{n-3}}\biggl(\sum_{j=1}^{\alpha_{n-2}} x_{n-2}^{\alpha_{n-2}-j}\delta_{n-2}(\sigma_{n-2}^{j-1}(\sigma_{n-1}^{\alpha_{n-1}}(\sigma_n^{\alpha_n}(r))))x_{n-2}^{j-1}\biggr)x_{n-1}^{\alpha_{n-1}}x_n^{\alpha_n}\\
+ &\ \dotsb + x_1^{\alpha_1}\biggl( \sum_{j=1}^{\alpha_2}x_2^{\alpha_2-j}\delta_2(\sigma_2^{j-1}(\sigma_3^{\alpha_3}(\sigma_4^{\alpha_4}(\dotsb (\sigma_n^{\alpha_n}(r))))))x_2^{j-1}\biggr)x_3^{\alpha_3}x_4^{\alpha_4}\dotsb x_{n-1}^{\alpha_{n-1}}x_n^{\alpha_n} \\
+ &\ \sigma_1^{\alpha_1}(\sigma_2^{\alpha_2}(\dotsb (\sigma_n^{\alpha_n}(r))))x_1^{\alpha_1}\dotsb x_n^{\alpha_n}, \ \ \ \ \ \ \ \ \ \ \sigma_j^{0}:={\rm id}_R\ \ {\rm for}\ \ 1\le j\le n.
\end{align*}
\end{proposition}
%%%%%%%%%%%%%%%
\begin{remark}[\cite{Reyes2015}, Remark 2.10)]\label{juradpr}
About Proposition \ref{coefficientes}, we have the following observation: Using (i), it follows that for the product $a_iX_ib_jY_j$, if $X_i:=x_1^{\alpha_{i1}}\dotsb x_n^{\alpha_{in}}$ and $Y_j:=x_1^{\beta_{j1}}\dotsb x_n^{\beta_{jn}}$, then when we compute every summand of $a_iX_ib_jY_j$ we obtain pro\-ducts of the coefficient $a_i$ with several evaluations of $b_j$ in $\sigma$'s and $\delta$'s depending of the coordinates of $\alpha_i$. This assertion follows from the expression:
%%%%%%%%%%%%%%%%
\begin{align*}
a_iX_ib_jY_j = &\ a_i\sigma^{\alpha_{i}}(b_j)x^{\alpha_i}x^{\beta_j} + a_ip_{\alpha_{i1}, \sigma_{i2}^{\alpha_{i2}}(\dotsb (\sigma_{in}^{\alpha_{in}}(b)))} x_2^{\alpha_{i2}}\dotsb x_n^{\alpha_{in}}x^{\beta_j} \\
+ &\ a_i x_1^{\alpha_{i1}}p_{\alpha_{i2}, \sigma_3^{\alpha_{i3}}(\dotsb (\sigma_{{in}}^{\alpha_{in}}(b)))} x_3^{\alpha_{i3}}\dotsb x_n^{\alpha_{in}}x^{\beta_j} \\
+ &\ a_i x_1^{\alpha_{i1}}x_2^{\alpha_{i2}}p_{\alpha_{i3}, \sigma_{i4}^{\alpha_{i4}} (\dotsb (\sigma_{in}^{\alpha_{in}}(b)))} x_4^{\alpha_{i4}}\dotsb x_n^{\alpha_{in}}x^{\beta_j}\\
+ &\ \dotsb + a_i x_1^{\alpha_{i1}}x_2^{\alpha_{i2}} \dotsb x_{i(n-2)}^{\alpha_{i(n-2)}}p_{\alpha_{i(n-1)}, \sigma_{in}^{\alpha_{in}}(b)}x_n^{\alpha_{in}}x^{\beta_j} \\
+ &\ a_i x_1^{\alpha_{i1}}\dotsb x_{i(n-1)}^{\alpha_{i(n-1)}}p_{\alpha_{in}, b}x^{\beta_j}.
\end{align*}
\end{remark}
%%%%%%%%%%%%%%%%%%%%
%%%%%%%%%%%%%%%%%%%%
Next, we recall the notion of $(\Sigma, \Delta)$-compatibility for rings.
%%%%%%%%%%%%%%%
\begin{definition}[\cite{ReyesSuarezUMA2017}, Definition 3.2]\label{Definition3.52008}
Consider a ring $R$ with a family of endomorphisms  $\Sigma$ and a family of $\Sigma$-derivations $\Delta$. Then,
\begin{enumerate}
\item $R$ is said to be $\Sigma$-{\em compatible}, if for each $a,b\in R$, $a\sigma^{\alpha}(b)=0$ if and only if $ab=0$, for every $\alpha\in \mathbb{N}^{n}$;
\item $R$ is said to be $\Delta$-{\em compatible}, if for each $a,b \in R$,  $ab=0$ implies $a\delta^{\beta}(b)=0$, for every $\beta \in \mathbb{N}^{n}$.
\end{enumerate}
If $R$ is both $\Sigma$-compatible and $\Delta$-compatible, $R$ is called $(\Sigma, \Delta)$-{\em compatible}.
\end{definition}
%%%%%%%%%%%%%%%%
\begin{examples}
Next, we present remarkable examples of $\sigma$-PBW extensions over $(\Sigma, \Delta)$-com\-pa\-ti\-ble rings (see \cite{Reyes2013PhD}, \cite{LezamaReyes2014} or \cite{ReyesYesica} for a detailed definition and reference of every example).
\begin{enumerate}
\item [\rm (a)] If $A$ is a skew PBW extension of a ring $R$ where the coefficients commute with the variables, that is, $x_ir = rx_i$, for every $r\in R$ and each $i=1,\dotsc, n$, or equivalently, $\sigma_i = {\rm id}_R$ and $\delta_i = 0$, for every $i$ (these extensions were called {\em constant} by the author in \cite{SuarezReyes2017}),  then it is clear that $R$ is $(\Sigma,\Delta)$-compatible. Some examples of constant $\sigma$-PBW extensions are the following: PBW extensions defined by Bell and Goodearl (which include the classical commutative polynomial rings, universal enveloping algebra of a Lie algebra, and others); some operator algebras (for example, the algebra of linear partial differential operators, the algebra of linear partial shift operators, the algebra of linear partial difference operators, the algebra of linear partial $q$-dilation operators, and the algebra of linear partial q-differential operators); the class of di\-ffu\-sion algebras; Weyl algebras; additive analogue of the Weyl algebra; multiplicative analogue of the Weyl algebra; some quantum Weyl algebras as $A_2(J_{a,b})$; the quantum algebra $\cU'(\mathfrak{so}(3,\Bbbk))$; the family of 3-dimensional skew polynomial algebras (there are exactly fifteen of these algebras, see \cite{ReyesSuarez2017FEJM}); Dispin algebra $\cU(osp(1,2))$; Woronowicz algebra $\cW_v(\mathfrak{sl}(2,\Bbbk))$; the complex algebra $V_q(\mathfrak{sl}_3(\mathbb{C}))$; $q$-Heisenberg algebra ${\bf H}_n(q)$; the Hayashi algebra $W_q(J)$, and more.
\item [\rm (b)] We also encounter examples of $\sigma$-PBW extensions (which are not constant) over $(\Sigma, \Delta)$-com\-pa\-ti\-ble rings. Let us see: (i) the quantum plane $\cO_q(\Bbbk^{2})$; the algebra of $q$-differential operators $D_{q,h}[x,y]$; the mixed algebra $D_h$; the operator differential rings; the algebra of differential operators $D_{\bf q}(S_{\bf q})$ on a quantum space ${S_{\bf q}}$, and more.
\item [\rm (c)] It is important to say that several algebras of quantum physics can be expressed as skew PBW extensions (for instance, Weyl algebras, additive and multiplicative analogue of the Weyl algebra, quantum Weyl algebras, $q$-Heisenberg algebra, and others), which allows us to characterize several properties with physical meaning. As Curado et. al., say in \cite{Curadoetal2008},  \textquotedblleft algebraic methods have long been applied to the solution of a large number of quantum physical systems. In the last decades, quantum algebras appeared in the framework of quantum integrable one-dimensional models and have ever since been applied to many physical phenomena [...] It was found that it could be generalized leading to the concept of deformed Heisenberg algebras that have been used in many areas, as nuclear physics, condensed matter, atomic physics, etc\textquotedblright. With these ideas in mind, next, we present some remarkable examples of these algebras (the proof that these algebras are skew PBW extensions can be realized using the theory developed in \cite{ReyesSuarezMomento2017}) which are $(\Sigma,\Delta)$-compatibles.
\begin{enumerate}
\item [\rm (i)] The Lie-deformed Heisenberg algebra introduced by Jannussis in \cite{Jannussisetal1992} is defined by the commutation relations
{\normalsize{\begin{align}
q_j(1+i\lambda_{jk})p_k - p_k(1-i\lambda_{jk})q_j = &\ i\hslash \delta_{jk}\notag \\
[q_j, q_k] = &\ [p_j, p_k] = 0,\ \ j, k = 1, 2, 3,\label{perritu}
\end{align}}}
where $q_j, p_j$ are the position and momentum operators, and $\lambda_{jk} = \lambda_{k} \delta_{jk}$, with $\lambda_k$ real parameters. If $\lambda_{jk} = 0$ one recovers the usual Heisenberg algebra.
\item [\rm (ii)] The {\em quantum Weyl algebra} introduced by Giaquinto and Zhang in \cite{GiaquintoZhang1995} with the aim of studying the Jordan Hecke symmetry is as a quantization of the usual second Weyl algebra. By definition, $A_2(J_{a,b})$ is the $\Bbbk$-algebra generated by the variables
$x_1,x_2,\partial_1,\partial_2$, with relations (depending on
parameters $a,b\in \Bbbk$)
{\normalsize{\begin{align}
x_1x_2 & = x_2x_1+ax_1^2, \ \ \ \ \ \ \ \ \ \ \ \ \ \ \ \partial_2\partial_1 = \partial_1\partial_2+b\partial_2^2\notag \\
\partial_1x_1 & = 1+x_1\partial_1 +ax_1\partial_2, \ \ \ \ \ \ \ \partial_1x_2 = -ax_1\partial_1-abx_1\partial_2+x_2\partial_1+bx_2\partial_2\notag \\
\partial_2x_1 & =x_1\partial_2, \ \ \ \ \ \ \ \ \ \ \ \ \ \ \ \ \ \ \ \ \ \ \ \ \partial_2x_2 = 1-bx_1\partial_2+x_2\partial_2.\label{aaaaaaa}
\end{align}}}
Over any field $\Bbbk$, if $a=b=0$, then $A_2(J_{0,0})\cong A_2$, the usual second Weyl algebra.
\item [\rm (iii)] With the purpose of obtaining bosonic representations of the Drinfield-Jimbo quantum algebras, Hayashi considered in \cite{Hayashi1990} the algebra {\bf U}. Let us see its construction (we follow \cite{Berger1992}, Example 2.7.7). Let {\bf U} be the algebra generated by the indeterminates $\omega_1,\dotsc, \omega_n, \psi_1,\dotsc, \psi_n,\psi_1^{*},\dotsc, \psi_n^{*}$, with the relations
{\small{\begin{align}
\psi_j \psi_i - \psi_i \psi_j =  &\ \psi_j^{*}\psi_i^{*} - \psi_i^{*}\psi_{j}^{*} = \omega_j\omega_i - \omega_i\omega_j = \psi_j^{*}\psi_{i} -  \psi_{i}\psi_j^{*} = 0, &\ 1\le i < j\le n,\notag \\
\omega_j\psi_i - q^{-\delta_{ij}}\psi_i\omega_j = &\ \psi_{j}^{*}\omega_i - q^{-\delta_{ij}}\omega_i\psi_j^{*} = 0, &\ 1\le i, j\le n,\notag \\
\psi_i^{*}\psi_{i} - q^{2}\psi_{i}\psi_{i}^{*} = &\
-q^{2}\omega_i^{2}, \ \ \ q\in \mathbb{C} &\ 1\le i \le n.
\label{claroo}
\end{align}}}
\item [\rm (iv)] Jannussis et. al., \cite{Jannussisetal1995} studied the non-Hermitian realization of a Lie deformed, a non-canonical Heisenberg algebra, considering the case of operators $A_j,\ B_k$ which are non-Hermitian (i.e., $\hslash=1$)
{\normalsize{\begin{align}
A_j(1+i\lambda_{jk})B_k - B_k(1-i\lambda_{jk})A_j = &\ i\delta_{jk}\notag \\
[A_j, B_k] = &\ 0\ \ (j\neq k)\notag \\
[A_j, A_k] = &\ [B_j, B_k] = 0,\label{perrituu}
\end{align}}}
and,
{\normalsize{\begin{align}
A_j^{+}(1 + i\lambda_{jk})B_k^{+} - B_k^{+}(1-i\lambda_{jk})A_j^{+} = &\ i\delta_{jk}\notag \\
[A_j^{+}, B_k^{+}] = &\ 0\ (j\neq k),\notag \\
[A_j^{+}, A_k^{+}] = &\ [B_j^{+}, B_k^{+}] = 0\label{perrituur}
\end{align}}}
where $A_j\neq A_j^{+},\ B_k \neq B_k^{+}$ $(j, k = 1, 2, 3)$. If the operators $A_j,\ B_k$ are in the form $A_j = f_j(N_j + 1)a_j,\ B_k = a_k^{+}f_k(N_k+1)$, where $a_j,\ a_j^{+}$ are leader operators of the usual Heisenberg-Weyl algebra, with $N_j$ the corresponding number operator $(N_j = a_j^{+}a_j,\ N_j\mid n_j\rangle = n_j|n_j\rangle)$, and the structure functions $f_j(N_j+1)$ complex, then it is showed in \cite{Jannussisetal1995} that $A_j$ and $B_k$ are given by
{\normalsize{\begin{align*}
A_j = &\ \sqrt{\frac{i}{1+i\lambda_j}} \biggl(\frac{[(1-i\lambda_j)/(1+i\lambda_j)]^{N_{j}+1}-1}{(1-i\lambda_j)/(1+i\lambda_j)-1} \frac{1}{N_j+1}\biggr)^{\frac{1}{2}}a_j\\
B_k = &\ \sqrt{\frac{i}{1+i\lambda_k}}a_k^{+} \biggl(\frac{[(1-i\lambda_k)/(1+i\lambda_k)]^{N_k + 1}-1}{(1-i\lambda_k)/(1+i\lambda_k)-1}\frac{1}{N_k+1}\biggr)^{\frac{1}{2}}.
\end{align*}}}
\end{enumerate}
\end{enumerate}
\end{examples}
%%%%%%%%%%%%%
With the aim of establishing the key results of the paper, we recall the following pre\-li\-mi\-na\-ry results  about $(\Sigma, \Delta)$-compatible rings.
%%%%%%%%%%%%%%%%%
\begin{proposition}[\cite{ReyesSuarezUMA2017}, Proposition 3.7]\label{colosss}
Let $R$ be an $(\Sigma, \Delta)$-compatible ring. For every $a, b \in R$, we have the following assertions:
\begin{enumerate}
\item [\rm (i)] if $ab=0$, then $a\sigma^{\theta}(b) = \sigma^{\theta}(a)b=0$, for each $\theta\in \mathbb{N}^{n}$.
\item [\rm (ii)] If $\sigma^{\beta}(a)b=0$ for some $\beta\in \mathbb{N}^{n}$, then $ab=0$.
\item [\rm (iii)] If $ab=0$, then $\sigma^{\theta}(a)\delta^{\beta}(b)= \delta^{\beta}(a)\sigma^{\theta}(b) = 0$, for every $\theta, \beta\in \mathbb{N}^{n}$.
\end{enumerate}
\end{proposition}
%%%%%%%%%%%%%%%%%%%%%%%
The next theorem generalizes \cite{Ouyang2010}, Proposition 3.6. We need to assume that the elements $c_{i,j}$ of Definition \ref{gpbwextension} (iv) are central in $R$.
%%%%%%%%%%%%
\begin{theorem}\label{misererer}
If $A$ is a $\sigma$-PBW extension of a reversible and $(\Sigma, \Delta)$-compatible ring $R$, then for every element $f=\sum_{i=0}^{m} a_iX_i\in A$, $f\in {\rm nil}(A)$ if and only if $a_i\in {\rm nil}(R)$, for each $1\le i\le m$.
\begin{proof}
Let $f\in A$ given as above and suppose that $f\in {\rm nil}(A)$ (consider $X_1 \prec X_2 \prec \dotsb \prec X_m$). Consider the notation established in Proposition \ref{coefficientes}. There exists a positive integer $k$ such that $f^{k} = (a_0 + a_1X_1 + \dotsb + a_mX_m)^{k}=0$. As an illustration, note that
\begin{align*}
f^{2} = &\ (a_mX_m + \dotsb + a_1X_1 + a_0)(a_mX_m + \dotsb + a_1X_1 + a_0)\\
= &\ a_mX_ma_mX_m +\ {\rm other\ terms\ less\ than}\ {\rm exp}(X_m)\\
= &\ a_m[\sigma^{\alpha_m}(a_m)X_m + p_{\alpha_m, a_m}]X_m + \ {\rm other\ terms\ less\ than}\ {\rm exp}(X_m)\\
= &\ a_m\sigma^{\alpha_m}(a_m)X_mX_m + a_mp_{\alpha_m, a_m}X_m +\ {\rm other\ terms\ less\ than}\ {\rm exp}(X_m)\\
= &\ a_m\sigma^{\alpha_m}(a_m)[c_{\alpha_m, \alpha_m}x^{2\alpha_m} + p_{\alpha_m, \alpha_m}] + a_mp_{\alpha_m, a_m}X_m +\ {\rm other\ terms\ less\ than}\ {\rm exp}(X_m)\\
= &\ a_m\sigma^{\alpha_m}(a_m)c_{\alpha_m, \alpha_m}x^{2\alpha_m} + \ {\rm other\ terms\ less\ than}\ {\rm exp}(x^{2\alpha_m}),
\end{align*}
and hence,
\begin{align*}
f^{3} = &\ (a_m\sigma^{\alpha_m}(a_m)c_{\alpha_m, \alpha_m}x^{2\alpha_m} + \ {\rm other\ terms\ less\ than}\ {\rm exp}(x^{2\alpha_m})) (a_mX_m + \dotsb + a_1x_1 + a_0)\\
= &\ a_m\sigma^{\alpha_m}(a_m)c_{\alpha_m, \alpha_m}x^{2\alpha_m}a_mX_m + \ {\rm other\ terms\ less\ than}\ {\rm exp}(x^{3\alpha_m})\\
= &\ a_m\sigma^{\alpha_m}(a_m)c_{\alpha_m, \alpha_m}[\sigma^{2\alpha_m}(a_m)x^{2\alpha_m} + p_{2\alpha_m, a_m}]X_m + \ {\rm other\ terms\ less\ than}\ {\rm exp}(x^{3\alpha_m})\\
= &\ a_m\sigma^{\alpha_m}(a_m)c_{\alpha_m,\alpha_m}\sigma^{2\alpha_m}(a_m)x^{2\alpha_m}X_m + \ {\rm other\ terms\ less\ than}\ {\rm exp}(x^{3\alpha_m})\\
= &\ a_m\sigma^{\alpha_m}(a_m)c_{\alpha_m,\alpha_m}\sigma^{2\alpha_m}(a_m)[c_{2\alpha_m, \alpha_m}x^{3\alpha_m} + p_{2\alpha_m,\alpha_m}]\\
= &\ a_m\sigma^{\alpha_m}(a_m)c_{\alpha_m,\alpha_m}\sigma^{2\alpha_m}(a_m)c_{2\alpha_m, \alpha_m}x^{3\alpha_m} + \ {\rm other\ terms\ less\ than}\ {\rm exp}(x^{3\alpha_m}).
\end{align*}
Continuing in this way, one can show that for $f^{k}$,
\[
f^{k} = a_m\prod_{l=1}^{k-1}\sigma^{l\alpha_m}(a_m)c_{l\alpha_m, \alpha_m}x^{k\alpha_m} + \ {\rm other\ terms\ less\ than}\ {\rm exp}(x^{k\alpha_m}),
\]
whence $0={\rm lc}(f^{k}) = a_m\prod_{l=1}^{k-1}\sigma^{l\alpha_m}(a_m)c_{l\alpha_m, \alpha_m}$, and since the elements $c$'s are central in $R$ and left invertible (Proposition \ref{coefficientes}), we have $0={\rm lc}(f^{k}) = a_m\prod_{l=1}^{k-1}\sigma^{l\alpha_m}(a_m)$. Using the $\Sigma$-compatibility of $R$, we obtain $a_m\in {\rm nil}(R)$.

Now, since
\begin{align*}
f^{k} = &\ ((a_0 + a_1X_1 + \dotsb + a_{m-1}X_{m-1}) + a_mX_m)^{k} \\
= &\ ((a_0 + a_1X_1 + \dotsb + a_{m-1}X_{m-1}) + a_mX_m)((a_0 + a_1X_1 + \dotsb + a_{m-1}X_{m-1}) + a_mX_m)\\
&\ \dotsb ((a_0 + a_1X_1 + \dotsb + a_{m-1}X_{m-1}) + a_mX_m)\ \ \ \ \ \ (k\ {\rm times})\\
= &\ [(a_0 + a_1X_1 + \dotsb + a_{m-1}X_{m-1})^{2} + (a_0 + a_1X_1 + \dotsb + a_{m-1}X_{m-1})a_mX_m\\
&\ + a_mX_m(a_0 + a_1X_1 + \dotsb + a_{m-1}X_{m-1}) + a_mX_ma_mX_m]\\
&\ \dotsb ((a_0 + a_1X_1 + \dotsb + a_{m-1}X_{m-1}) + a_mX_m)\\
= &\ (a_0 + a_1X_1 + \dotsb + a_{m-1}X_{m-1})^{k} + h,
\end{align*}
where $h$ is an element of $A$ which involves products of monomials with the term $a_mX_m$ on the left and the right, by  Proposition \ref{lindass}, Remark \ref{juradpr} and having in mind that $a_m\in {\rm nil}(R)$, which is an ideal of $R$ (remember that reversible implies semicommutative), the expression for $f^{k}$ reduces to $f^k = (a_0 + a_1X_1 + \dotsb + a_{m-1}X_{m-1})^{k}$. Using a similar reasoning as above, one can prove that
\[
f^{k} = a_{m-1}\prod_{l=1}^{k-1}\sigma^{l(\alpha_{m-1})}(a_{m-1})c_{l(\alpha_{m-1}), \alpha_{m-1}}x^{k\alpha_{m-1}} + \ {\rm other\ terms\ less\ than}\ {\rm exp}(x^{k\alpha_{m-1}}).
\]
Hence ${\rm lc}(f^{k}) = a_{m-1}\prod_{l=1}^{k-1}\sigma^{l\alpha_{m-1}}(a_{m-1})c_{l\alpha_{m-1}, \alpha_{m-1}}$, and so $a_{m-1}\in {\rm nil}(R)$. If we repeat this argument, it follows that $a_i\in {\rm nil}(R)$, for $0\le i\le m$.

Conversely, suppose that $a_i\in {\rm nil}(R)$, for every $i$. If $k_i$ is the minimum integer positive such that $a_i^{k_i} = 0$, for every $i$, let $k:={\rm max}\{k_i\mid 1\le i\le n\}$. It is clear that $a_i^{k}=0$, for all $i$. Let us prove that $f^{(m+1){k}+1} = 0$, and hence, $f\in {\rm nil}(A)$. Since the  expression for $f$ have $m+1$ terms, when we realize the product $f^{(m+1){k}+1}$ we have sums of products of the form
\begin{equation}\label{rigoo}
a_{i,1}X_{i,1}a_{i,2}X_{i,2}\dotsb a_{i, (m+1){k}}X_{i, (m+1){k}}a_{i,(m+1){k}+1}X_{i,(m+1){k}+1}.
\end{equation}
Note that there are exactly $(m+1)^{(m+1)k+1}$ products of the form (\ref{rigoo}). Now, since when we compute $f^{(m+1){k}+1}$ every product as (\ref{rigoo}) involves at least $k$ elements $a_i$, for some $i$, then every one of these products is equal to zero by Proposition \ref{lindass}, Remark \ref{juradpr} and the $(\Sigma, \Delta)$-compatibility of $R$ (more exactly, Proposition \ref{colosss}). In this way, every term of $f^{(m+1){k}+1}$ is equal to zero, and hence $f\in {\rm nil}(A)$.
\end{proof}
\end{theorem}
%%%%%%%%%%%%%%%%%%%%%%%
\section{Skew $\Pi$-Armendariz rings}\label{reunionMelendez}
In \cite{LunqunJinwangyYueming2013}, Definition 2.1, Lunqun et. al., introduced the notion of skew $\pi$-Armendariz ring in the following way: let $B$ be a ring with an endomorphism $\alpha$ and an $\alpha$-derivation $\delta$. $B$ is called a {\em skew} $\pi$-{\em Armendariz ring}, if for polynomials $f(x) = \sum_{i=0}^{m} a_ix^{i}$ and $g(x) = \sum_{j=0}^{n} b_jx_j$ in $B[x;\alpha,\delta]$, $f(x)g(x)\in {\rm nil}(B[x;\alpha,\delta])$ implies $a_ib_j\in {\rm nil}(B)$, for each $0\le i \le m$ and $0\le j\le n$. Of course, every $\alpha$-Armendariz ring defined by Hong et. al., \cite{HongKwakRizvi2006} is skew $\pi$-Armendariz considering $\delta$ as the zero derivation. If $B$ is a ring with an endomorphism $\alpha$ and an $\alpha$-derivation $\delta$, following Moussavi and Hashemi \cite{MoussaviHashemi2005}, $B$ is said to be a $(\alpha,\delta)$-{\em skew Armendariz ring}, if for polynomials $f(x) = \sum_{i=0}^{m} a_ix_i$ and $g(x)=\sum_{j=0}^{t}$ in the Ore extension $B[x;\alpha,\delta]$, $f(x)g(x)=0$ implies $a_ix^{i}b_jx^{j}=0$, for each $i,\ j$. In  \cite{LunqunJinwangyYueming2013}, Theorem 2.6, it was proved that skew $\pi$-Armendariz rings are more general than skew Armendariz rings when the ring of coefficients $B$ is $(\alpha, \delta)$-compatible.

With the aim of extending the above definition for the context of $\sigma$-PBW extensions, next we present the following notion:
%%%%%%%%%%%%%%%
\begin{definition}\label{skewpiArmendariz}
Let $A$ be a $\sigma$-PBW extension of a ring $R$. $R$ is called a {\em skew}-$\Pi$ {\em Armendariz ring}, if for elements $f=\sum_{i=0}^{m} a_iX_i,\ g=\sum_{j=0}^{t} b_jY_j$ of $A$, $fg\in {\rm nil}(A)$ implies that $a_ib_j\in {\rm nil}(R)$, for every $0\le i\le m$ and $0\le j\le t$.
\end{definition}
%%%%%%%%%%%%%
Now, since in \cite{NinoReyes2017}, Definition 3.4 (see Definition \ref{ninoreyesdef2017} below), the author introduced the notion of $(\Sigma,\Delta)$-{\em skew Armendariz} over $\sigma$-PBW extensions (this definition  generalizes the treatments developed for both classical polynomial rings and Ore extensions of injective type (for example \cite{Armendariz1974}, \cite{RegeChhawchharia1997}, \cite{AndersonCamillo1998}, \cite{HongKimKwak2003}, \cite{Matczuk2004}), the natural purpose is to establish the version of \cite{LunqunJinwangyYueming2013}, Theorem 2.6 for the general case of $\sigma$-PBW extensions over $(\Sigma,\Delta)$-skew Armendariz rings. This result is formulated in Proposition  \ref{frtet}.
%Nasr-Isfahani et. al., \cite{NasrMoussavi2008}, $B$ is said to be {\em skew Armendariz}, if for polynomials $f(x) = \sum_{i=0}^{m} a_ix^{i},\ g(x) = \sum_{j=0}^{n} b_jx^{j}\in B[x;\alpha,\delta]$, $f(x)g(x)$ implies $a_0b_j=0$, for each $0\le j \le n$.
%%%%%%%%%%%%%%%
\begin{definition}[\cite{NinoReyes2017}, Definition 3.4]\label{ninoreyesdef2017}
If $A$ is a $\sigma$-PBW extension of a ring $R$, then $R$ is called a $(\Sigma, \Delta)$-{\em skew Armendariz ring}, if whenever $f = \sum_{i=0}^{t} a_iX_i,\ g = \sum_{j=0}^{s} b_jY_j\in A$ with $fg=0$, then $a_iX_ib_jY_j=0$, for every value of $i$ and $j$.
\end{definition}

The next proposition establishes that $(\Sigma,\Delta)$-skew Armendariz which are $\Sigma$-compatible are skew $\Pi$-Armendariz. We assume that the elements $c_{i,j}$ of Definition \ref{gpbwextension} (iv) are central in $R$.
%%%%%%%%%%%%%%%%%%
\begin{proposition}\label{frtet}
Let $A$ be a $\sigma$-PBW extension of a ring $R$. If $R$ is $\Sigma$-compatible and $(\Sigma, \Delta)$-{\em skew Armendariz ring}, then $R$ is skew $\Pi$-Armendariz.
\begin{proof}
First of all, let us prove that if we have $p_1, \dotsc, p_l$ elements of $A$ with $p_1\dotsb p_l = 0$, then if $a_k\in C_{p_k}$, for $k=1,\dotsc, l$, we have $a_1\dotsb a_l=0$. We proceed by induction following the notation considered in Proposition \ref{coefficientes}. If $l=2$, let $p_1=\sum_{i=0}^{m} a_iX_i,\ p_2=\sum_{j=0}^{t} b_jY_j$. By assumption we have $a_iX_ib_jY_j = 0$, for every $i, j$, and since $a_iX_ib_jY_j = a_i[\sigma^{\alpha_i}(b_j)X_i + p_{\alpha_i, b_j}]Y_j = a_i\sigma^{\alpha_i}(b_j)X_iY_j + a_ip_{\alpha_i, b_j}Y_j = a_i\sigma^{\alpha_i}(b_j)[c_{\alpha_i,\beta_j}x^{\alpha_i + \beta_j} + p_{\alpha_i, \beta_j}] + a_ip_{\alpha_i, b_j}Y_j = a_i\sigma^{\alpha_i}(b_j)c_{\alpha_i, \beta_j}x^{\alpha_i+\beta_j} + a_i\sigma^{\alpha_i}(b_j)p_{\alpha_i,\beta_j} + a_ip_{\alpha_i, b_j}Y_j$, we obtain ${\rm lc}(a_iX_ib_jY_j) = a_i\sigma^{\alpha_i}(b_j) = 0$, for each $i, j$ (using that the elements $c_{\alpha_i, \beta_j}$ are central in $R$). By the $\Sigma$-compatibility of $R$, $a_ib_j = 0$, for every value of $i$ and $j$, so the assertion follows.

Second of all, let $l > 2$. If $h:=p_2p_3\dotsb p_l$, then $p_1h=0$, and by the reasoning above, $a_1a_h=0$, where $a_1\in C_{p_1},\ a_h\in C_{h}$. Having in mind the form of the elements of $h$, that is, $a_h =a_2\dotsb a_l$, where $a_2\in C_{f_2}, \dotsc, a_l\in C_{f_l}$ (which is due to the fact that $R$ is $(\Sigma, \Delta)$-skew Armendariz and $\Sigma$-compatible), then we obtain $a_1\dotsb a_l=0$.

Finally, if we have two elements $f, g\in A$ given by $f = \sum_{i=0}^{t} a_iX_i$ and $g = \sum_{j=0}^{s} b_jY_j$ with $fg\in {\rm nil}(A)$, then there exists a positive integer $r$ with $(fg)^{r} = 0$ whence by the analysis above, $a_ib_j\in {\rm nil}(R)$, for $i=1,\dotsc, t,\ j=1,\dotsc, s$, and hence $R$ is skew $\Pi$-Armendariz.
\end{proof}
\end{proposition}
%%%%%%%%%%%%%%%
A more general notion than $(\Sigma, \Delta)$-skew Armendariz ring it was formulated in \cite{ReyesSuarez2016}, Definition 3.1. More exactly,
%%%%%%%%%%%%%
\begin{definition}[\cite{ReyesSuarez2016b}, Definition 3.1]
Let $A$ be a skew PBW extension of a ring $R$. $R$ is called a $\Sigma$-{\em skew Armendariz ring}, if for elements $f=\sum_{i=0}^{m} a_iX_i$ and $g=\sum_{j=0}^{t} b_jY_j$ in $A$, the equality $fg=0$ implies $a_i\sigma^{\alpha_i}(b_j) = 0$, for all $0\le i\le m$ and $0\le j\le t$, where $\alpha_i = {\rm exp}(X_i)$.
\end{definition}
%%%%%%%%%%%%%
We can extend the result established in Proposition \ref{frtet} to $\Sigma$-skew Armendariz rings (again, we assume that the elements $c_{i,j}$ of Definition \ref{gpbwextension} (iv) are central in $R$):
%%%%%%%%%%%%%%%%
\begin{proposition}\label{1234}
Let $A$ be a $\sigma$-PBW extension of a ring $R$ with a family of endomorphisms $\Sigma = \{\sigma_1,\dotsc, \sigma_n\}$ and a family of $\Sigma$-derivations $\Delta = \{\delta_1,\dotsc, \delta_n\}$. If $R$ is $\Sigma$-compatible and $
\Sigma$-skew Armendariz ring, then $R$ is skew $\Pi$-Armendariz.
\begin{proof}
Let us prove that if we have $p_1, \dotsc, p_l$ elements of $A$ with $p_1\dotsb p_l = 0$, then if $a_k\in C_{p_k}$, for $k=1,\dotsc, l$, we have $a_1\dotsb a_l=0$. We proceed by induction. If $l=2$, let $p_1=\sum_{i=0}^{m} a_iX_i,\ p_2=\sum_{j=0}^{t} b_jY_j$. By assumption we have $a_i\sigma^{\alpha_i}(b_j)=0$, for every $i, j$, and by the $\Sigma$-compatibility of $R$, $a_ib_j = 0$, for every value of $i$ and $j$. The rest of the proof is completely similar to the presented in Proposition \ref{frtet}.
\end{proof}
\end{proposition}
%%%%%%%%%%%%%%
\begin{remark}
Recently, in \cite{ReyesSuarezClifford2017}, Definition 4.1, it was introduced the notion of {\em skew Armendariz} ring which is more general than both $(\Sigma, \Delta)$-skew Armendariz and $\Sigma$-skew Armendariz. More precisely,  if $R$ is a ring and $A$ is a  $\sigma$-PBW extension of $R$, we say that $R$ is a {\em skew-Armendariz} ring, if for polynomials $f=a_0+a_1X_1+\dotsb + a_mX_m$ and $g=b_0+b_1Y_1 + \dotsb + b_tY_t$ in $A$, $fg=0$ implies $a_0b_k=0$, for each $0\le k\le t$. We believe that the following assertion is true: if $A$ is a $\sigma$-PBW extension of a ring $R$ with a family of endomorphisms $\Sigma = \{\sigma_1,\dotsc, \sigma_n\}$ and a family of $\Sigma$-derivations $\Delta = \{\delta_1,\dotsc, \delta_n\}$, where $R$ is $(\Sigma, \Delta)$-compatible and  skew Armendariz, then $R$ is skew $\Pi$-Armendariz. In a forthcoming paper we will investigate this conjecture.
\end{remark}
%%%%%%%%%%%%%%%%%
Next, we present some results concerning NI rings and its relation with skew $\Pi$-Ar\-men\-da\-riz rings. We denote ${\rm nil}(R)A:=\{f\in A\mid f= a_0 + a_1X_1 + \dotsb + a_mX_m,\ a_i\in {\rm nil}(R)\}$. We start with the following two useful results which can be considered as the analogous results to \cite{Ouyang2010}, Lemma 3.4  and Lemma 3.5, respectively.
%%%%%%%%%%%%%%%
\begin{lemma}\label{pajaritocabas}
If $ab\in {\rm nil}(R)$, where $R$ is a $(\Sigma,\Delta)$-compatible and reversible ring, then $a\sigma^{\alpha}(\delta^{\beta}(b))$ and $a\delta^{\beta}(\sigma^{\alpha}(b))$ are elements of ${\rm nil}(R)$.
% In this way, if $A$ is a $\sigma$-PBW extension of $A$ and $X, Y\in {\rm Mon}(A)$, then $aXbY$ is an element of ${\rm nil}(R)A$.
\begin{proof}
By assumption there exists a positive integer $k$ such that $(ab)^{k}=0$. Consider the following equalities:
\begin{align*}
(ab)^{k} = &\ ab ab \dotsb ab ab ab\ \ \ (k\ {\rm times})\\
= &\ abab \dotsb ab ab a\sigma^{\alpha}(\delta^{\beta}(b))\ \ \ ({\rm Proposition\ \ref{colosss}\ (iii)})\\
= &\ a\sigma^{\alpha}(\delta^{\beta}(b)) ab ab ab \dotsb ab ab\ \ \ (R\ {\rm is\ reversible})\\
= &\ a\sigma^{\alpha}(\delta^{\beta}(b)) ab ab \dotsb ab a\sigma^{\alpha}(\delta^{\beta}(b))\ \ \ ({\rm Proposition\ \ref{colosss}\ (iii)})\\
= &\ a\sigma^{\alpha}(\delta^{\beta}(b)) a\sigma^{\alpha}(\delta^{\beta}(b)) ab ab \dotsb ab\ \ \ (R\ {\rm is\ reversible})\\
\end{align*}
Following this procedure we guarantee that the element $a\sigma^{\alpha}(\delta^{\beta}(b))$ belongs to ${\rm nil}(R)$. For the element $a\delta^{\beta}(\sigma^{\alpha}(b))$ the reasoning is completely similar. %Finally, note that $aXbY$
\end{proof}
\end{lemma}
%%%%%%%%%%%
\begin{remark}
In \cite{LunqunJinwangyYueming2013}, Lemma 2.3 (2) there is a little mistake. There it was omitted the condition of reversibility of the ring as we can appreciate in the original reference \cite{Ouyang2010}, Lemma 3.4.
\end{remark}
%%%%%%%%%%%%
\begin{lemma}\label{beatles}
If $R$ is a $(\Sigma,\Delta)$-compatible ring, then $a\sigma^{\theta}(b)\in {\rm nil}(R)$ implies $ab\in {\rm nil}(R)$, for every $\theta \in \mathbb{N}^{n}$.
\begin{proof}
Since $a\sigma^{\theta}(b)\in {\rm nil}(R)$, there exists a positive integer $k$ with $(a\sigma^{\theta}(b))^{k}=0$. We have the following assertions
\begin{align*}
(a\sigma^{\theta}(b))^{k} =  &\ a\sigma^{\theta}(b) a\sigma^{\theta}(b) \dotsb a\sigma^{\theta}(b) a\sigma^{\theta}(b)\ \ \ (k\ {\rm times})\\
= &\ a\sigma^{\theta}(b)a\sigma^{\theta}(b) \dotsb a\sigma^{\theta}(b) ab \ \ \ ({\rm Definition\ of}\ \Sigma-{\rm compatibility})\\
= &\ a\sigma^{\theta}(b) a\sigma^{\theta}(b) \dotsb a\sigma^{\theta}(b) \sigma^{\theta}(ab)\ \ \ ({\rm Proposition\ \ref{colosss}\ (i)})\\
= &\ a\sigma^{\theta}(b) a\sigma^{\theta}(b) \dotsb a \sigma^{\theta}(bab)\ \ \ (\sigma^{\theta}\ {\rm is\ an\ endomorphism\ of}\ R)\\
= &\ a\sigma^{\theta}(b) a\sigma^{\theta}(b) \dotsb abab\ \ \
({\rm Definition\ of}\ \Sigma-{\rm compatibility})
\end{align*}
If we continue in this way, we can see that the element $ab\in {\rm nil}(R)$, which concludes the proof.
\end{proof}
\end{lemma}
%%%%%%%%%%%%%
\begin{theorem}\label{hijoamado}
Let $A$ be a skew PBW extension of $R$. If $R$ is $(\Sigma, \Delta)$-compatible and reversible, then $R$ is skew $\Pi$-Armendariz.
\begin{proof}
First of all, let us show that if $R$ is a $(\Sigma,\Delta)$-compatible NI ring, then ${\rm nil}(A)\subseteq {\rm nil}(R)A$. Consider an element $f\in {\rm nil}(A)$ given by the expression $f= a_0 + a_1X_1 + \dotsb + a_mX_m$, with $X_1\prec X_2\prec \dotsb \prec X_m$. Then, there exists a positive integer $k$ with $f^{k} = 0$. As an illustration, note that
\begin{align*}
f^{2} = &\ (a_mX_m + \dotsb + a_1x_1 + a_0)(a_mX_m + \dotsb + a_1x_1 + a_0)\\
= &\ a_mX_ma_mX_m +\ {\rm other\ terms\ less\ than}\ {\rm exp}(X_m)\\
= &\ a_m[\sigma^{\alpha_m}(a_m)X_m + p_{\alpha_m, a_m}]X_m + \ {\rm other\ terms\ less\ than}\ {\rm exp}(X_m)\\
= &\ a_m\sigma^{\alpha_m}(a_m)X_mX_m + a_mp_{\alpha_m, a_m}X_m +\ {\rm other\ terms\ less\ than}\ {\rm exp}(X_m)\\
= &\ a_m\sigma^{\alpha_m}(a_m)[c_{\alpha_m, \alpha_m}x^{2\alpha_m} + p_{\alpha_m, \alpha_m}] + a_mp_{\alpha_m, a_m}X_m +\ {\rm other\ terms\ less\ than}\ {\rm exp}(X_m)\\
= &\ a_m\sigma^{\alpha_m}(a_m)c_{\alpha_m, \alpha_m}x^{2\alpha_m} + \ {\rm other\ terms\ less\ than}\ {\rm exp}(x^{2\alpha_m}),
\end{align*}
and hence,
\begin{align*}
f^{3} = &\ (a_m\sigma^{\alpha_m}(a_m)c_{\alpha_m, \alpha_m}x^{2\alpha_m} + \ {\rm other\ terms\ less\ than}\ {\rm exp}(x^{2\alpha_m})) (a_mX_m + \dotsb + a_1x_1 + a_0)\\
= &\ a_m\sigma^{\alpha_m}(a_m)c_{\alpha_m, \alpha_m}x^{2\alpha_m}a_mX_m + \ {\rm other\ terms\ less\ than}\ {\rm exp}(x^{3\alpha_m})\\
= &\ a_m\sigma^{\alpha_m}(a_m)c_{\alpha_m, \alpha_m}[\sigma^{2\alpha_m}(a_m)x^{2\alpha_m} + p_{2\alpha_m, a_m}]X_m + \ {\rm other\ terms\ less\ than}\ {\rm exp}(x^{3\alpha_m})\\
= &\ a_m\sigma^{\alpha_m}(a_m)c_{\alpha_m,\alpha_m}\sigma^{2\alpha_m}(a_m)x^{2\alpha_m}X_m + \ {\rm other\ terms\ less\ than}\ {\rm exp}(x^{3\alpha_m})\\
= &\ a_m\sigma^{\alpha_m}(a_m)c_{\alpha_m,\alpha_m}\sigma^{2\alpha_m}(a_m)[c_{2\alpha_m, \alpha_m}x^{3\alpha_m} + p_{2\alpha_m,\alpha_m}]\\
= &\ a_m\sigma^{\alpha_m}(a_m)c_{\alpha_m,\alpha_m}\sigma^{2\alpha_m}(a_m)c_{2\alpha_m, \alpha_m}x^{3\alpha_m} + \ {\rm other\ terms\ less\ than}\ {\rm exp}(x^{3\alpha_m}).
\end{align*}
Continuing in this way, one can show that for $f^{k}$,
\[
f^{k} = a_m\prod_{l=1}^{k-1}\sigma^{l\alpha_m}(a_m)c_{l\alpha_m, \alpha_m}x^{k\alpha_m} + \ {\rm other\ terms\ less\ than}\ {\rm exp}(x^{k\alpha_m}),
\]
whence $0={\rm lc}(f^{k}) = a_m\prod_{l=1}^{k-1}\sigma^{l\alpha_m}(a_m)c_{l\alpha_m, \alpha_m}$, and since the elements $c$'s are central in $R$ and left invertible (Proposition \ref{coefficientes}), we have $0={\rm lc}(f^{k}) = a_m\prod_{l=1}^{k-1}\sigma^{l\alpha_m}(a_m)$. Using the $\Sigma$-compatibility of $R$, we obtain $a_m\in {\rm nil}(R)$.

Now, since
\begin{align*}
f^{k} = &\ ((a_0 + a_1X_1 + \dotsb + a_{m-1}X_{m-1}) + a_mX_m)^{k} \\
= &\ ((a_0 + a_1X_1 + \dotsb + a_{m-1}X_{m-1}) + a_mX_m)((a_0 + a_1X_1 + \dotsb + a_{m-1}X_{m-1}) + a_mX_m)\\
&\ \dotsb ((a_0 + a_1X_1 + \dotsb + a_{m-1}X_{m-1}) + a_mX_m)\ \ \ \ k\ {\rm times}\\
= &\ [(a_0 + a_1X_1 + \dotsb + a_{m-1}X_{m-1})^{2} + (a_0 + a_1X_1 + \dotsb + a_{m-1}X_{m-1})a_mX_m\\
&\ + a_mX_m(a_0 + a_1X_1 + \dotsb + a_{m-1}X_{m-1}) + a_mX_ma_mX_m]\\
&\ \dotsb ((a_0 + a_1X_1 + \dotsb + a_{m-1}X_{m-1}) + a_mX_m)\\
= &\ (a_0 + a_1X_1 + \dotsb + a_{m-1}X_{m-1})^{k} + h,
\end{align*}
where $h$ is an element of $A$ which involves products of monomials with the term $a_mX_m$ on the left and the right. From Proposition \ref{lindass}, Remark \ref{juradpr} and having in mind that $a_m\in {\rm nil}(R)$ which is an ideal of $R$ ($R$ is reversible and hence a NI ring), the expression for $f^{k}$ reduces to $f^k = (a_0 + a_1X_1 + \dotsb + a_{m-1}X_{m-1})^{k}$, and using a similar reasoning as above, one can prove that
\[
f^{k} = a_{m-1}\prod_{l=1}^{k-1}\sigma^{l(\alpha_{m-1})}(a_{m-1})c_{l(\alpha_{m-1}), \alpha_{m-1}}x^{k\alpha_{m-1}} + \ {\rm other\ terms\ less\ than}\ {\rm exp}(x^{k\alpha_{m-1}}).
\]
Hence ${\rm lc}(f^{k}) = a_{m-1}\prod_{l=1}^{k-1}\sigma^{l\alpha_{m-1}}(a_{m-1})c_{l\alpha_{m-1}, \alpha_{m-1}}$, and so $a_{m-1}\in {\rm nil}(R)$. If we repeat this argument, then $a_i\in {\rm nil}(R)$, for $0\le i\le m$.

Finally, let us prove that $R$ is skew $\Pi$-Armendariz. Consider two elements $f, g \in A$ given by $f=\sum_{i=0}^{m} a_iX_i$ and $g=\sum_{j=0}^{t} b_jY_j$ with $fg\in {\rm nil}(A)$. Note that
\begin{equation*}
fg = \sum_{k=0}^{m+t} \biggl( \sum_{i+j=k} a_iX_ib_jY_j\biggr) \in {\rm nil}(A)\subseteq {\rm nil}(R)A,
\end{equation*}
and ${\rm lc}(fg)= a_m\sigma^{\alpha_m}(b_t)c_{\alpha_m, \beta_t}\in {\rm nil}(R)$. Since the elements $c_{i,j}$ are in the center of $R$, then $c_{\alpha_m,\beta_t}$ are also in the center of $R$, whence $a_m\sigma^{\alpha_m}(b_t)\in {\rm nil}(R)$, and by Lemma \ref{beatles} it follows that $a_mb_t\in {\rm nil}(R)$. The idea is to prove that $a_pb_q\in {\rm nil}(R)$, for $p+q\ge 0$.  We proceed  by induction. Suppose that $a_pb_q\in {\rm nil}(R)$, for $p+q=m+t, m+t-1, m+t-2, \dotsc, k+1$, for some $k>0$. By Lemma \ref{pajaritocabas},  we obtain $a_pX_pb_qY_q\in {\rm nil}(R)A$ for these values of $p+q$. In this way, it is sufficient to consider the sum of the products $a_uX_ub_vY_v$, where $u+v=k, k-1,k-2,\dotsc, 0$. Fix $u$ and $v$. Consider the sum of all terms  of $fg$  having exponent $\alpha_u+\beta_v$. By Proposition \ref{lindass}, Remark \ref{juradpr} and the assumption $fg\in {\rm nil}(A)$, we know that the sum of all coefficients of all these terms  can be written as
\begin{equation}\label{Federer}
a_u\sigma^{\alpha_u}(b_v)c_{\alpha_u, \beta_v} + \sum_{\alpha_{u'} + \beta_{v'} = \alpha_u + \beta_v} a_{u'}\sigma^{\alpha_{u'}} ({\rm \sigma's\ and\ \delta's\ evaluated\ in}\ b_{v'})c_{\alpha_{u'}, \beta_{v'}}\in {\rm nil}(R).
\end{equation}
As we suppose above, $a_pb_q\in {\rm nil}(R)$ for $p+q=m+t, m+t-1, \dotsc, k+1$, so Lemma \ref{pajaritocabas} guarantees that the product $a_p(\sigma'$s and $\delta'$s evaluated in $b_q$), for any order of $\sigma'$s and $\delta'$s, is an element of ${\rm nil}(R)$. Since $R$ is reversible, then $({\rm \sigma's\ and\ \delta's\ evaluated\ in}\ b_{q})a_p\in {\rm nil}(R)$. In this way, multiplying (\ref{Federer}) on the right by $a_k$, and using the fact that the elements $c$'s are in the center of $R$, we obtain
\begin{equation}\label{doooooo}
a_u\sigma^{\alpha_u}(b_v)a_kc_{\alpha_u, \beta_v} + \sum_{\alpha_{u'} + \beta_{v'} = \alpha_u + \beta_v} a_{u'}\sigma^{\alpha_{u'}} ({\rm \sigma's\ and\ \delta's\ evaluated\ in}\ b_{v'})a_kc_{\alpha_{u'}, \beta_{v'}}\in {\rm nil}(R),
\end{equation}
whence, $a_u\sigma^{\alpha_u}(b_0)a_k\in {\rm nil}(R)$. Since $u+v=k$ and $v=0$, then $u=k$, so $a_k\sigma^{\alpha_k}(b_0)a_k\in {\rm nil}(R)$, from which $a_k\sigma^{\alpha_k}(b_0)\in {\rm nil}(R)$ and hence $a_kb_0=0$ by Lemma \ref{beatles}. Therefore, we now have to study the expression (\ref{Federer}) for $0\le u \le k-1$ and $u+v=k$. If we multiply (\ref{doooooo}) on the right by $a_{k-1}$, then
{\small{\[
a_u\sigma^{\alpha_u}(b_v)a_{k-1}c_{\alpha_u, \beta_v} + \sum_{\alpha_{u'} + \beta_{v'} = \alpha_u + \beta_v} a_{u'}\sigma^{\alpha_{u'}} ({\rm \sigma's\ and\ \delta's\ evaluated\ in}\ b_{v'})a_{k-1}c_{\alpha_{u'}, \beta_{v'}}\in {\rm nil}(R).
\]}}
Using a similar reasoning as above, we can see that $a_u\sigma^{\alpha_u}(b_1)a_{k-1}c_{\alpha_u, \beta_1}\in {\rm nil}(R)$. Since the elements $c$'s are central and left invertible, $a_u\sigma^{\alpha_u}(b_1)a_{k-1}\in {\rm nil}(R)$, and using the fact $u=k-1$, we have $a_{k-1}\sigma^{\alpha_{k-1}}(b_1)\in {\rm nil}(R)$, from which $a_{k-1}b_1\in {\rm nil}(R)$. Continuing in this way we prove that $a_ib_j\in {\rm nil}(R)$, for $i+j=k$. Therefore $a_ib_j\in {\rm nil}(R)$, for $0\le i\le m$ and $0\le j\le t$.
\end{proof}
\end{theorem}
%%%%%%%%%%%%
\begin{corollary}\label{sport}
If $A$ is a $\sigma$-PBW extension of $R$, then every $\Sigma$-rigid ring is skew $\Pi$-Armendariz.
\begin{proof}
From \cite{Reyes2015}, p. 182, we know that $\Sigma$-rigid rings are reduced, so reversible. Now, by \cite{ReyesSuarezUMA2017}, Proposition 3.3, $\Sigma$-rigid rings are $(\Sigma, \Delta)$-compatible rings. The assertion follows from Theorem \ref{hijoamado}.
\end{proof}
\end{corollary}
%%%%%%%%%%%%%%%
%Finally, in \cite{AndersonCamillo1998},  Theorem 2, Anderson and Camillo proved that a ring $B$ is Armendariz if and only if the classical commutative ring $B[x]$ is Armendariz. For the rings of interest in this paper, we have the following result which extends \cite{LunqunJinwangyYueming2013}, Theorem 3.7.
%%%%%%%%%%%%%%%%%
%\begin{theorem}
%Let $A$ be a $\sigma$-PBW extension of a ring $R$. If $R$ is a $(\Sigma, \Delta)$-compatible 2-primal ring, then $R[x]$ is a skew $\pi$-Armendariz ring.
%\end{theorem}
%%%%%%%%%%%
\begin{remark}\label{masabiert}
In \cite{RodriguezReyes2017}, Theorem 3.6, under the same conditions in Theorem \ref{hijoamado}, the author proved that $R$ is a $(\Sigma,\Delta)$-skew McCoy ring (see \cite{RodriguezReyes2017} for the definition of these rings). This fact shows the importance of assuming both conditions ($(\Sigma, \Delta)$-compatibility and reversibility) on $R$. Now, having in mind the results obtained in this paper (Propositions \ref{frtet} and \ref{1234}, and Corollary \ref{sport}), it is interesting to investigate the properties of Baer, quasi-Baer, p.p. and p.q.-Baer, considering the notion of skew $\Pi$-Armendariz with the aim of extending the theory developed in \cite{Reyes2015}, \cite{NinoReyes2017},  \cite{ReyesSuarez2016b} and \cite{ReyesYesica}.
\end{remark}
%%%%%%%%%%%%%%%%%
%%%%%%%%%%%%%%%%%
%%%%%%%%%%%%%%%%%%%%%%%%%%%%%%%

\vspace{0.5cm}

\noindent {\bf \Large{Acknowledgements}}

\vspace{0.5cm}

The author was supported by the research fund of Facultad de Ciencias, Universidad Nacional de Colombia, Bogot\'a, Colombia, HERMES CODE 30366.

%%%%%%%%%%%%%%%%%%%%%%%%%%%%%%%%%%%%%%%%%%%%%%%%%%%%%%%%%%%%%%%%%%

\end{document}